\documentclass{amsart}
\usepackage{amssymb}
\usepackage{amsmath}

\usepackage{a4}

\begin{document}
\def\id{\operatorname{id}}
\def\bird#1{[#1]}
\def\mapright#1{\smash{\mathop{\longrightarrow}\limits\sp{#1}}}
\newtheorem{theorem}{Theorem}[section]
\newtheorem{lemma}[theorem]{Lemma}
\newtheorem{remark}[theorem]{Remark}
\newtheorem{definition}[theorem]{Definition}
\newtheorem{corollary}[theorem]{Corollary}
\newtheorem{example}[theorem]{Example}
\def\qedbox{\hbox{$\rlap{$\sqcap$}\sqcup$}}
\makeatletter
  \renewcommand{\theequation}{%
   \thesection.\alph{equation}}
  \@addtoreset{equation}{section}
 \makeatother
\font\pbglie=eufm10
\def\gg{\text{\pbglie{g}}}
\title[Lie Groups]
{The spectral geometry of the canonical Riemannian submersion of a compact Lie Group}
\author{C. Dunn, P. Gilkey, and J.H. Park${}^1$}

\begin{address}{CD: Mathematics Department, California State University at San Bernardino,
San Bernardino, CA 92407, USA. Email: \it
cmdunn@csusb.edu.}\end{address}
\begin{address}{PG: Mathematics Department, University of Oregon, Eugene, OR 97403, USA.
Email:
\it gilkey@uoregon.edu.}\end{address}
\begin{address}{JP: Department of Mathematics, SungKyunKwan University,
Suwon, 440-746, SOUTH KOREA. E-mail: \it parkj@skku.edu}\end{address}
\begin{abstract}Let $G$ be a compact connected Lie group which is equipped with a
bi-invariant Riemannian metric. Let $m(x,y)=xy$
be the multiplication operator. We show the associated fibration $m:G\times G\rightarrow G$
is a Riemannian
submersion with totally geodesic fibers and we study the spectral geometry of this
submersion. We show the pull back of eigenforms on the base have finite Fourier series on the total space  and we give
examples where arbitrarily many Fourier coefficients can be non-zero.  We give necessary and sufficient conditions
that the pull back of a form on the base is harmonic on the total space.
\end{abstract}
\keywords{ Riemannian submersion, eigenform, finite fourier
series
\newline 1991 {\it Mathematics Subject Classification.} 58G25\newline ${}^1$ Corresponding author}
\maketitle

\section{Introduction} The spectral geometry of Riemannian manifolds has been studied
extensively; compact
Lie groups play a central role in this investigation. For example, work of  Schueth
\cite{Sc01} shows there are
non-trivial isospectral families of left invariant metrics on compact Lie groups, although
any such family which
includes a bi-invariant metric is necessarily trivial; this work has been extended by Proctor \cite{P05}.
Riemannian submersions of Lie groups  with totally geodesic
fibers have been studied by Ranjan \cite{R86}.
We refer to
\cite{GLP99} for a further  discussion of the
spectral geometry of Riemannian submersions.

In this paper, we shall study the spectral geometry of the
multiplication map $m:G\times G\rightarrow G$ where $G$ is a
compact Lie group. Adopt the following notational conventions. Let
$\Delta_M^p:=(d\delta+\delta d)$ be the Laplace-Beltrami operator
acting on the space of smooth $p$ forms $C^\infty(\Lambda^p(M))$
on a compact smooth closed Riemannian manifold $M$ of dimension
$m$. We summarize briefly the following well known facts which we
shall need, see, for example \cite{G04} for further details.
Denote the distinct eigenvalues and associated eigenspaces by:
\begin{eqnarray*}
&&\operatorname{Spec}(\Delta_M^p)=\{0=\lambda_0<\lambda_1<...<\lambda_n<...\},\\
&&E_\lambda(\Delta_M^p)=\left\{\phi\in
C^\infty(\Lambda^p(M)):\Delta_M^p\phi=\lambda\phi\right\}\,.
\end{eqnarray*}
The spectral multiplicities $\dim\{E_\lambda(\Delta_M^p)\}$ are all finite. Furthermore there is a complete
orthonormal decomposition
$$
L^2(\Lambda^p(M))=\oplus_{\lambda\in\operatorname{Spec}(\Delta_M^p)}
E_\lambda(\Delta_M^p)\,.
$$

Let $G$ be a compact connected Lie group which is equipped with a bi-invariant Riemannian
metric
$ds^2_G$. Normalize the product metric on $G\times G$ by taking
\begin{equation}\label{eqn-1.a}
ds^2_{G\times G}=2(ds^2_G\oplus ds^2_G)\,.
\end{equation}

The situation on $0$-forms is particularly simple; we shall show in Section \ref{sect-2}
that the pull-back of an
eigenfunction is again an eigenfunction with the same eigenvalue:

\begin{theorem}\label{thm-1.1}
Let $ds^2_G$ be a bi-invariant metric on a compact Lie group $G$. Let
$ds^2_{G\times G}=2(ds^2_G\oplus ds^2_G)$. Then the multiplication map
$m:G\times G\rightarrow G$ is a Riemannian submersion with totally geodesic fibers and
$m^*\{E_\lambda(\Delta_G^0)\}\subset E_\lambda(\Delta_{G\times G}^0)$.
\end{theorem}

Let $\pi_\lambda$ be orthogonal projection on $E_\lambda(\Delta_M^p)$. If
$\phi\in C^\infty(\Lambda^p(M))$, let $\mu(\phi)$ be the number of eigenvalues $\lambda$ so
that $\pi_\lambda\phi\ne0$;
this is the number of distinct eigenvalues which are involved in the Fourier series
decomposition of $\phi$. We shall
use the Peter-Weyl theorem in Section \ref{sect-3} to show that:

\begin{theorem}\label{thm-1.2}
Let $ds^2_G$ be a bi-invariant metric on a compact Lie group $G$. Let
$ds^2_{G\times G}=2(ds^2_G\oplus ds^2_G)$. If $\phi\in E_\lambda(\Delta_G^p)$, then
$\mu(m^*\phi)\le\textstyle{{2\dim\{G\}}\choose{p}}\dim\{E_\lambda(\Delta_G^p)\}$.
\end{theorem}

The geometry of left invariant $1$-forms plays a central role in our discussions. The
following result will be
established in Section \ref{sect-4}:

\begin{theorem}\label{thm-1.3}
Let $ds^2_G$ be a bi-invariant metric on a compact Lie group $G$. Let
$ds^2_{G\times G}=2(ds^2_G\oplus ds^2_G)$. Let $\phi\in E_\lambda(\Delta_G^1)$ be left invariant. Then one
may decompose
$m^*\phi=\Phi_1+\Phi_2$ where $0\ne\Phi_1\in E_{\frac32\lambda}(\Delta_{G\times G}^1)$ and
$0\ne\Phi_2\in E_{\frac12\lambda}(\Delta_{G\times G}^1)$.
\end{theorem}

Theorem \ref{thm-1.2} shows that the pull back of an eigenform
has a finite Fourier  series. In Section \ref{sect-5},
we will use Theorem
\ref{thm-1.3} to establish following result which shows that the number of eigenvalues
involved in the Fourier
decomposition of
$m^*\phi$ can be arbitrarily large:

\begin{theorem}\label{thm-1.4}
Let $p\ge1$ and let $\mu_0\in\mathbb{N}$ be given. There exists a bi-invariant metric on a compact
Lie group
$G$,
there
exists $\lambda$, and there exists $0\ne\phi\in E_\lambda(\Delta_G^p)$ so that
$\mu(m^*\phi)=\mu_0$.\end{theorem}

The Hodge-DeRham theorem identifies the $n^{\operatorname{th}}$
cohomology group
$H^n(M;\mathbb{C})$ of $M$ with the space of harmonic $n$-forms $E_0(\Delta_M^n)$ if $M$ is a
compact Riemannian manifold. Thus the eigenvalue $0$ has a particular significance. Let
$\Lambda(E_0(\Delta_G^1))$ be the subring generated over $\mathbb{C}$
by the harmonic $1$-forms; one has that $\phi\in\Lambda^n(E_0(\Delta_G^1))$ if and only if one can
express:
$$
  \phi=\sum_{|I|=n}a_I\phi^{i_1}\wedge...\wedge\phi^{i_n}\text{ where }a_I\in\mathbb{C}\text{ and }
  \phi^{i}\in E_0(\Delta_G^1)\,.
$$
\begin{theorem}\label{thm-1.5}
Let $ds^2_G$ be a bi-invariant metric on a compact Lie group $G$. Let
$ds^2_{G\times G}=2(ds^2_G\oplus ds^2_G)$. Assume $G$ connected.
\begin{enumerate}
\item $\Lambda^n(E_0(\Delta_G^1))\subset E_0(\Delta_G^n)$.
\item $\phi\in\Lambda^n(E_0(\Delta_G^1))$ if and only if $m^*\phi\in E_0(\Delta_{G\times G}^n)$.
\item Let $G$ be simply connected. If $\phi\in
E_0(\Delta_G^n)$ for $n>0$,
$m^*\phi\notin E_0(\Delta_{G\times G}^n)$.
\end{enumerate}
\end{theorem}

One can consider more generally the situation where $G$ and $G\times G$ are
endowed with arbitrary left invariant metrics $ds^2_G$ and $ds^2_{G\times G}$ where there is no a priori
relation assumed between these metrics. The question of when this is a Riemannian submersion is an
interesting one and will be studied in more detail in a subsequent paper. For the moment, however,
we content ourselves in Section
\ref{sect-1.7} by generalizing Theorem
\ref{thm-1.2} to this setting:

\begin{theorem}\label{thm-1.6}
Let $G$ and $G\times G$ be equipped with left invariant metrics
$ds^2_G$ and
$ds^2_{G\times G}$. If $\phi\in E_\lambda(\Delta_G^p)$, then
$$
\mu(m^*\phi)\le
{\textstyle\binom{2\dim\{G\}}{p}}^2{\textstyle\binom{\dim\{G\}}{p}}^2\dim\{E_\lambda(\Delta_G^p)\}^4\,.
$$
\end{theorem}

We remark that this bound is much worse than the
bound given in Theorem \ref{thm-1.2}; at $2$ different points in the proof we shall need to pass from a left
invariant subspace to a biinvariant subspace and this greatly increases estimate on the dimension.

\section{The geometry of the multiplication map $m$}\label{sect-2}

Let $\pi:X\rightarrow Y$ be a surjective smooth map where $X$ and $Y$ are compact
Riemannian manifolds. We suppose
that $\pi$ is a submersion, i.e. that the map $\pi_*:T_xX\rightarrow T_{\pi x}Y$ is
surjective for every $x\in X$, and
let
$\mathcal{V}$ (resp. $\mathcal{H}$) be the associated vertical (resp. horizontal)
distribution:
$$
\mathcal{V}:=\{\xi\in TX:\pi_*\xi=0\}\quad\text{and}\quad
\mathcal{H}:=\mathcal{V}^\perp\,.
$$
We say that $\pi$ is a {\it Riemannian submersion} if $\pi_*:\mathcal{H}_x\rightarrow
TY_{\pi x}$ is an isometry
$\forall x$.

The following example is instructive. Let $m(u,v)=u+v$ define a linear map from
$\mathbb{R}^{2n}\rightarrow\mathbb{R}^n$. Take the standard Euclidean metric on
$\mathbb{R}^n$. We may identify
$T_x\mathbb{R}^{2n}=\mathbb{R}^{2n}$ and $T_y\mathbb{R}^n=\mathbb{R}^n$. Under this
identification,
$$\mathcal{V}=\textstyle\operatorname{Span}_{\xi\in\mathbb{R}^n}\{(\frac12\xi,-\frac12\xi)\}
\quad\text{and}\quad
  \mathcal{H}=\operatorname{Span}_{\xi\in\mathbb{R}^n}\{(\frac12\xi,\frac12\xi)\}\,.$$
We have $m_*(\frac12\xi,\frac12\xi)=\xi$. Thus if $\xi$ is a unit vector in
$T_y\mathbb{R}^{n}$, we need that
$(\frac12\xi,\frac12\xi)$ is a unit vector in $T_x\mathbb{R}^{2n}$. This motivates the
factor of $2$ which appears in
Equation (\ref{eqn-1.a}) since the ordinary Euclidean length of $(\frac12\xi,\frac12\xi)$
would be $\frac12$ and not
$1$. With this normalization, $m$ becomes a Riemannian submersion.

More generally, let $G$ be a Lie group which is equipped with a bi-invariant
Riemannian metric $ds^2_G$. Let $m(x,y)=xy$ be the multiplication operator from $G\times
G\rightarrow G$. Let
$\{e_i^L\}$ (resp. $\{e_i^R\}$) be an orthonormal frame of left (resp.
right) invariant vector fields on $G$. We assume $e_i^L(1)=e_i^R(1)=e_i$ where $1\in G$ is
the
unit of the group and where $\{e_i\}$ is an orthonormal basis for $T_1(G)$. Let $\exp$ be
the exponential map in the
group. Then the flows
$\Xi_i^L$ and
$\Xi_i^R$ of these vector fields are:
$$
\Xi_i^L:(g,t)\rightarrow g\exp
(te_i)\quad\text{and}\quad\Xi_i^R:(g,t)\rightarrow\exp(te_i)g\,.
$$

The multiplication map $m$ defines a smooth surjective map $m:G\times G\rightarrow G$.
Consider the
following curves in $G\times G$ with initial position $(g_1,g_2)$:
\begin{eqnarray*}
&&\gamma_i^{g_1,g_2}:t\rightarrow(g_1\exp(\textstyle\frac12(te_i),
\exp(-\textstyle\frac12te_i)g_2),\\
&&\varrho_i^{g_1,g_2}:t\rightarrow(g_1\exp(\textstyle\frac12te_i),\exp(\textstyle\frac12te_i
)g_2),\\
&&\tau_i^{g_1,g_2}:t\rightarrow(\exp(te_i)g_1,g_2)\,.
\end{eqnarray*}

We may identify $T(G\times G)=TG\oplus TG$. Because
$m\tau_i^{g_1,g_2}:t\rightarrow\exp(te_i)g_1g_2$,
$$m_*\{\dot\tau_i^{g_1,g_2}(0)\}=e_i^R(m(g_1,g_2))\,.$$
Consequently
$m_*$ is surjective so
$m$ is a submersion. As $m\gamma_i^{g_1,g_2}:t\rightarrow g_1g_2$ is independent of $t$, one has
$\dot\gamma_i=\textstyle\frac12(e_i^L,-e_i^R)\in\ker\{m_*\}$. It now follows that
\begin{equation}\label{eqn-2.a}
\begin{array}{l}
\mathcal{V}:=\ker\{m_*\}=\operatorname{Span}\{V_i:=\textstyle\frac12(e_i^L,-e_i^R)\},\\
\mathcal{H}:=\ker\{m_*\}^\perp=\operatorname{Span}\{H_i:=\textstyle\frac12(e_i^L,e_i^R)\}\,.
\vphantom{\vrule height 11pt}\end{array}\end{equation}
Let $L_g$ and $R_g$ denote left and right  multiplication in the group. As
$\dot\varrho_i=H_i$,
\begin{equation}\label{eqn-2.b}
m_{*(g_1,g_2)} \{H_i(g_1,g_2)\}=(L_{g_1})_*(R_{g_2})_*e_i\,.
\end{equation}

Since $L_{g_1}$ and $R_{g_2}$ are isometries, it follows that $\{m_*H_i(g_1,g_2)\}$ is an
orthonormal basis for
$T_{g_1g_2}G$. We have defined $ds_{G\times G}^2=2(ds_{G}^2\oplus ds_G^2)$. We show that
$m$ is a Riemannian
submersion by computing:
\begin{eqnarray*}
&&(H_i,H_j)_{G\times G}=2\textstyle\frac14\{(e_i^L,e_j^L)_G+(e_i^R,e_j^R)_G\}=\delta_{ij},\\
&&(H_i,V_j)_{G\times G}=2\textstyle\frac14\{(e_i^L,e_j^L)_G-(e_i^R,e_j^R)_G\}=0,\\
&&(V_i,V_j)_{G\times G}=2\textstyle\frac14\{(e_i^L,e_j^L)_G+(e_i^R,e_j^R)_G\}=\delta_{ij}\,.
\end{eqnarray*}

Fix $h\in G$. The map
$T_h:(g_1,g_2)\rightarrow(hg_2^{-1},g_1^{-1}h)$ is an isometry $G\times G$. Clearly
$T_h(g_1,g_2)=(g_1,g_2)$ if
and only if $g_1=hg_2^{-1}$ and $g_2=g_1^{-1}h$ or equivalently if $g_1g_2=h$. Thus the
fixed point set of $T$ is
$m^{-1}(h)$. Since the fixed point set
of an isometry consists of the disjoint union of totally geodesic submanifolds, the fibers
of $m$, which are
connected submanifolds diffeomorphic to $G$, are totally geodesic.
It now follows that the mean curvature covector vanishes. Theorem 4.3.1 of \cite{GLP99}
shows
$m^*\Delta_G^0=\Delta_{G\times G}^0m^*$. This completes the
proof of Theorem \ref{thm-1.1}.

\section{The Peter-Weyl theorem}\label{sect-3}

We recall the classical Peter-Weyl theorem; for further details see, for example, \cite{H03,HM98}.
Let
$G$ be a compact Lie group which is equipped with a bi-invariant metric; assume the metric is
normalized so
$G$ has unit volume. If
$\rho$ is a smooth left representation of
$G$ on a finite dimensional complex vector space $V$, then by averaging an arbitrary innerproduct on $V$ over the
group we can always choose an innerproduct on $V$ which is preserved by $\rho$. Thus
any such representation is unitarizable. Let $\operatorname{Irr}(G)$ be the set of isomorphism
classes of finite dimensional irreducible unitary left representations of $G$.  We can decompose
any finite dimensional left representation space $V$ as a direct sum of irreducibles:
$$V=\oplus_{\rho\in\operatorname{Irr}(G)}n_\rho V_\rho\,;$$ the
multiplicities $n_\rho$ are independent of the particular decomposition chosen and are non-zero
for only finitely many $\rho$.

Let $\{e_i\}$ be an orthonormal basis for $V_\rho$ where $\rho\in\operatorname{Irr}(G)$. We may
expand $\rho(g)e_i=\sum_j\rho_{ij}(g)e_j$; the functions
$\rho_{ij}\in C^\infty(G)$ are said to be the matrix coefficients of $\rho$. We let
$$H_\rho:=\operatorname{Span}_{1\le i,j\le\dim(\rho)}\{\rho_{ij}\}\subset L^2(G)\,.$$
It is easily verified that $H_\rho$ is invariant under both the left and right group action and
that $H_\rho$ is independent of the particular orthonormal basis chosen for
$V_\rho$; isomorphic representations determine the same space. Furthermore, as a left
representation space for $G$, $H_\rho$ is isomorphic to $\dim(\rho)$ copies of the original
representation $\rho$.

If $V$ is any finite dimensional subspace of
$L^2(G)$ which is left-invariant under $G$ and which is abstractly isomorphic to $V_\rho$ as a
representation space, then one has $V\subset H_\rho$; to put it another way,
$H_\rho$ contains all the left submodules of
$L^2(G)$ which are isomorphic to $V_\rho$. Furthermore, we have a complete orthogonal direct sum decomposition
$$L^2(G)=\oplus_{\rho\in\operatorname{Irr}(G)}H_\rho=\oplus_{\rho\in\operatorname{Irr}(G)}\dim(\rho)\cdot
V_\rho\,.$$
This means that $\{\rho_{ij}\}_{1\le i,j\le\dim(\rho),\rho\in\operatorname{Irr}(G)}$ is
a complete orthonormal basis for $L^2(G)$.

More generally, let $\{\phi^i_L\}$ be an orthonormal basis for the space of left invariant $1$-forms. If one has that
$I=\{1\le i_1<...<i_p\le\dim(G)\}$ is a multi-index, let $\Phi^I_L:=\phi^{i_1}_L\wedge...\wedge\phi^{i_p}_L$;
the $\Phi^I_L$ are an orthonormal basis for the space of left invariant $p$-forms and as a left
representation space for $G$ one has:
\begin{equation}\label{eqn-3.a}
L^2(\Lambda^p(G))=\oplus_{\rho\in\operatorname{Irr}(G),|I|=p}H_\rho\otimes\Phi_L^I
=\oplus_{\rho\in\operatorname{Irr}(G)}{\textstyle\binom{\dim\{G\}}{p}}\dim(V_\rho)V_\rho\,.
\end{equation}
The subspace $H_\rho^p:=\oplus_{|I|=p}H_\rho\cdot\Phi_L^I$ is a bi-invariant $G$ submodule of
$L^2(\Lambda^p(G))$ which contain every left subrepresentation of $G$ on $L^2(\Lambda^p(G))$ isomorphic
to
$V_\rho$.

Let $\pi_\lambda$ be orthogonal projection from $L^2(\Lambda^p(G))$ to
$E_\lambda(\Delta_G^p)$ and let $\mu(\phi)$ be the number of eigenvalues $\lambda$ so
$\pi_\lambda(\phi)\ne0$. We
prepare for the proof of Theorem
\ref{thm-1.2} by establishing:

\begin{lemma}\label{lem-3.1}
Let $H\subset L^2(\Lambda^p(G))$ be invariant under the action of $L_g$ for all $g\in G$. If
$\phi\in H$, then $\mu(\phi)\le{\binom{\dim\{G\}}{p}}\dim\{H\}$.
\end{lemma}

\medbreak\noindent{\it Proof.}
Clearly $\pi_\lambda H$ is
non-trivial if and only if there exists
$\rho\in
\operatorname{Irr}(G)$ so that the multiplicities satisfy:
$$n_H(\rho)>0\quad\text{and}\quad n_{E_\lambda(\Delta_G^p)}(\rho)>0\,.$$
Note that only a finite number of representations appear in $H$ and
only a finite number of
eigenspaces involve any given representation. By Equation (\ref{eqn-3.a}),
\begin{eqnarray*}
&&\mu(\phi)\le
\sum_{\rho\in\operatorname{Irr}(G):n_\rho(H)\ne0}\bigg
\{\sum_{\lambda:n_\rho(E_\lambda(\Delta_G^p))\ne0}1\bigg\}\\
&\le&\sum_{\rho\in\operatorname{Irr}(G):n_\rho(H)\ne0}\bigg\{{{\dim\{G\}}\choose{p}}
\dim\{V_\rho\}\bigg\}
\le{{\dim\{G\}}\choose{p}}\dim\{H\}\,.\qquad\qedbox
\end{eqnarray*}

We can now establish Theorem \ref{thm-1.2}.  It is
convenient to introduce $\tilde m(g,h)=gh^{-1}$. Let
$H=E_\lambda(\Delta_G^p)$. Since the metric is bi-invariant, the
Laplacian and hence the eigenspaces are preserved by both left and
right multiplication. Let $\tilde H:=\tilde m^*H$. We compute:
\begin{eqnarray*}
&&\tilde m\{L_{g,h}^{G\times G}(a_1,a_2)\}=\tilde m(ga_1,ha_2)
   =ga_1a_2^{-1}h^{-1}=L_g^GR_{h^{-1}}^G\tilde m(a_1,a_2),\\
&&\{L_{g,h}^{G\times G}\}^*\tilde m^*=\tilde
m^*(R^G_{{h^{-1}}})^*(L_g^G)^*\,.
\end{eqnarray*}
Since $H$ is invariant under both the left and right actions of $G$, $\tilde H$ is invariant
under the left action of $G\times G$. We replace the group in question by $G\times G$ and apply Lemma \ref{lem-3.1}
to estimate
$\mu(\tilde m^*\phi)$. Since the
metric on
$G\times G$ is bi-invariant,
$\psi(x,y):=(x,y^{-1})$ is an isometry of
$G\times G$. We have
$$a_1a_2=m(a_1,a_2)=a_1(a_2^{-1})^{-1}=\tilde m(\psi(a_1,a_2))$$
and thus $m^*=\psi^*\tilde m^*$. Consequently $\mu(m^*\phi)=\mu(\psi^*\tilde
m^*\phi)=\mu(\tilde m^*\phi)$.
\hfill\qedbox

\section{Left invariant $1$-forms}\label{sect-4}
Let $\Lambda^p_L(G)$ be the finite dimensional vector space of left invariant $p$-forms on
$G$.
Define the left and right actions of $G$
on
$G\times G$ by:
\begin{equation}\label{eqn-4.a}
\begin{array}{ll}
L_{1,g}:(x,y)\rightarrow(gx,y),&L_{2,g}:(x,y)\rightarrow(x,gy),\\
R_{1,g}:(x,y)\rightarrow(xg,y),&R_{2,g}:(x,y)\rightarrow(x,yg)\vphantom{\vrule height 11pt}
\end{array}\end{equation}
Consider the following subspaces:
$$
\tilde\Lambda^p(G\times G)=\{\theta\in C^\infty(\Lambda^p(G\times G):
L_{1,g}^*\theta=\theta,R_{1,g^{-1}}^*L_{2,g}^*\theta=\theta\
\forall\ g\in G\}\,.
$$

\begin{lemma}\label{lem-4.1}
Adopt the notation established above. Then:
\begin{enumerate}
\item $d_{G\times G}\{\tilde\Lambda^p(G\times G)\}\subset\tilde\Lambda^{p+1}(G\times G)$,\qquad
      $\delta_{G\times G}\{\tilde\Lambda^{p+1}(G\times G)\}\subset\tilde\Lambda^{p}(G\times G)$,\newline
     $\Delta_{G\times G}^p\{\tilde\Lambda^p(G\times G)\}\subset\tilde\Lambda^p(G\times G)$, and
     $\tilde\Lambda^p(G\times G)\wedge\tilde\Lambda^q(G\times G)\subset\tilde\Lambda^{p+q}(G\times G)$.
\item The map $\theta\rightarrow\theta(1)$ is an isomorphism from $\tilde\Lambda^p(G\times G)$ to $\Lambda^p(G\times G)(1)$.
\item $m^*\{\Lambda_L^p(G)\}\subset \tilde\Lambda^p(G\times G)$.
\end{enumerate}
\end{lemma}

\begin{proof} Assertion (1) follows since
the maps of Equation (\ref{eqn-4.a}) are isometries and thus the pullbacks defined by these
maps commute with $d$, $\delta$,
$\Delta$, and $\wedge$. To prove Assertion (2), define an action $A$ of $G\times G$ on
$G\times G$ by setting:
$$A_{g,h}:(a,b)\rightarrow (gah^{-1},hb)$$ this is a fixed point free transitive isometric
group action since
$$A_{g_1,h_1}A_{g_2,h_2}=A_{g_1g_2,h_1h_2}\,.$$
This exhibits $G\times G$ as a homogeneous space.  We have furthermore that:
$$\begin{array}{lll}
m(ga,b)=gm(a,b),&m\circ L_{1,g}=L_g\circ m,&m^*L_g^*=L_{1,g}^*m^*,\\
m(ag^{-1},gb)=m(a,b),&m\circ L_{2,g} R_{1,g^{-1}}=m,&
R_{1,g^{-1}}^* L_{2,g}^*m^*=m^*\,.
\end{array}$$
Suppose that $\phi\in\Lambda^p_L(G)$. Then $L_g^*\phi=\phi$ for
all $g$. Consequently
$$L_{1,g}^*m^*\phi=m^*L_g^*\phi=m^*\phi\quad\text{and}\quad R_{1,g^{-1}}^*
L_{2,g}^*m^*\phi=m^*\phi\,.$$
Assertion (3) follows.\end{proof}

Fix an orthonormal frame $\{\phi_L^i\}$ for $\Lambda^1_L(G)$ so that
\begin{equation}\label{eqn-4.b}
\Delta_G^1\{\phi_L^i\}=\lambda_i\phi_L^i\,.
\end{equation}
Since right and left multiplication commute, right
multiplication preserves
$\Lambda_L^1(G)$. Thus we may decompose
\begin{equation}\label{eqn-4.c}
R_g^*\phi^i_L=\sum_j\xi_{ij}(g)\phi^j_L\,.
\end{equation}
Since $R_gR_h=R_{hg}$ and since $R_1=\operatorname{id}$, we have
$$\xi_{ij}(g)\xi_{jk}(h)=\xi_{ik}(hg)\quad\text{and}\quad\xi_{ij}(1)=\delta_{ij}\,.$$
We may decompose $\Lambda^1(G\times G)=\Lambda^1(G)\oplus\Lambda^1(G)$. Define
$$\Phi_1^i(u,v)=\sum_j\xi_{ij}(v)\phi_L^j(u)\oplus0\quad\text{and}\quad\Phi_2^i(u,v)=0\oplus
\phi_L^i(v)\,.$$

\begin{lemma}\label{lem-4.2}
Adopt the notation established above.
\begin{enumerate}
\item $\{\Phi_1^i,\Phi_2^i\}$ is a basis for $\tilde\Lambda^1(G\times G)$.
\item $m^*\phi^i_L=\Phi_1^i+\Phi_2^i$.
\item $\Delta_{G\times G}^1\Phi_1^i=\frac32\lambda_i\Phi_1^i$ and
$\Delta_{G\times G}^1\Phi_2^i=\frac12\lambda_i\Phi_2^i$.
\end{enumerate}
\end{lemma}

\begin{proof} It is immediate from the definition that $L_{1,g}^*\Phi_1^i=\Phi_1^i$,
$L_{1,g}^*\Phi_2^i=\Phi_2^i$, and $R_{1,g^{-1}}^*L_{2,g}^*\Phi_2^i=\Phi_2^i$. We use
Equation (\ref{eqn-4.c})
to see:
\begin{eqnarray*}
&&\{R_{1,g^{-1}}^*L_{2,g}^*\Phi_1^i\}(u,v)=\sum_{jk}\xi_{ij}(gv)\xi_{jk}(g^{-1})
\phi_L^k(u)\oplus0\\
&&\qquad=\sum_{jkl}\xi_{il}(v)\xi_{lj}(g)\xi_{jk}(g^{-1})\phi_L^k(u) \oplus0\\
&&\qquad=\sum_k\xi_{ik}(v)\phi_L^k(u)\oplus0=\Phi_1^i(u,v)\,.
\end{eqnarray*}
Thus $\Phi_1^i\in\tilde\Lambda^1(G\times G)$ and
$\Phi_2^i\in\tilde\Lambda^1(G\times G)$. Because $\Phi_1^i (1,
1)=\phi^i_L(1)\oplus0$ and because
$\Phi_2^i{(1,1)}=0\oplus\phi^i_L(1)$, Assertion (1) now follows
from Assertion (2) of Lemma \ref{lem-4.1}. We dualize Equations
(\ref{eqn-2.a}) and (\ref{eqn-2.b}) to see that
$$\{m^*\phi_L^i\}(1,1)=\phi_L^i(1)\oplus\phi_L^i(1)=\{\Phi_1^i+\Phi_2^i\}(1,1)\,.$$
The identity of Assertion (2) of Lemma \ref{lem-4.2} now follows
from Assertion (1) of Lemma \ref{lem-4.2} and from Assertion (3) of Lemma \ref{lem-4.1}.

Suppose $\phi\in\Lambda^1_L(G)$. Then $\delta_G\phi\in\Lambda^0_L(G)$ is left-invariant and
hence $\delta_G\phi=c$ is constant. Since
$dc=0$,

$$c^2\operatorname{vol}(G)=(\delta_G\phi,\delta_G\phi)_{L^2(G)}=(\phi,d_G\delta_G\phi)_
{L^2(\Lambda^1G)}=0\,.$$
Similarly if $\Phi\in\tilde\Lambda^1(G\times G)$, then $\delta_{G\times
G}\Phi\in\tilde\Lambda^0(G\times G)$ is
invariant under the transitive group action $A$ defined above. Consequently
$\delta_{G\times G}\Phi=C$ constant and
again
\begin{eqnarray*}
&&C^2\operatorname{vol}(G\times G)=(\delta_{G\times G}\Phi,\delta_{G\times
G}\Phi)_{L^2(G\times G)}
=(\Phi,d_{G\times G}\delta_{G\times G}\Phi)_{L^2\Lambda^1(G\times G)}=0\,.
\end{eqnarray*}
Consequently one may express:
\begin{equation}\label{eqn-4.d}
\Delta_G^1\{\phi_L^i\}=\delta_Gd_G\{\phi^i_L\}\quad\text{and}\quad
\Delta_{G\times G}^1\{\Phi^i_a\}=\delta_{G\times G}d_{G\times
G}\{\Phi^i_a\}\ \ \text{for}\ \ a=1,2\,.
\end{equation}

Decompose
$$\textstyle
d_G\{\phi_L^i\}=\sum_{j<k}C_{ijk}\phi_L^j\wedge\phi_L^k\quad\text{and}\quad
\delta_G\{\phi_L^j\wedge\phi_L^k\}=\sum_iD_{ijk}\phi_L^i\,.
$$
We compute:
\begin{eqnarray*}
&&D_{ijk}\text{vol}(G)=(\delta_G\{\phi_L^j\wedge\phi_L^k\},\phi_L^i)_{L^2(\Lambda^1G)}
=(\phi_L^j\wedge\phi_L^k,d\phi_L^i)_{L^2(\Lambda^2G)}\\
&=&C_{ijk}\text{vol}(G)\,.\end{eqnarray*}
Consequently $D_{ijk}=C_{ijk}$. Equations (\ref{eqn-4.b}) and (\ref{eqn-4.d}) yield:
$$
\sum_{j<k,l}C_{ljk}C_{ijk}\phi_L^l=\delta_G\bigg\{\sum_{j<k}C_{ijk}\phi_L^j\wedge\phi_L^k
\bigg\}
=\delta_Gd_G\{\phi_L^i\}=\Delta_G^1\{\phi_L^i\}
=\lambda_i\phi_L^i$$
and consequently
\begin{equation}\label{eqn-4.e}
\sum_{j<k}C_{ljk}C_{ijk}=\lambda_i\delta^{il}\,.
\end{equation}

Let $\sigma_2(g_1,g_2)=g_2$ denote projection on the second factor. Since $\Phi_2^i=\sigma_2^*\phi_L^i$ and since
$\Phi_1^i+\Phi_2^i=m^*\phi_L^i$,
\begin{equation}\label{eqn-4.f}
\begin{array}{l}
d_{G\times G}\{\Phi_2^i\}=\sum_{j<k}C_{ijk}\Phi_2^j\wedge\Phi_2^k,\\
d_{G\times
G}\{\Phi_1^i+\Phi_2^i\}=\sum_{j<k}C_{ijk}(\Phi_1^j+\Phi_2^j)\wedge(\Phi_1^k+\Phi_2^k),
\vphantom{\vrule height 12pt}\\
d_{G\times G}\{\Phi_1^i\}=d_{G\times G}\{\Phi_1^i+\Phi_2^i\}-d_{G\times
G}\{\Phi_2^i\}\vphantom{\vrule height 12pt}\\
  \qquad\qquad
=\sum_{j<k}C_{ijk}\bigg\{\Phi_1^j\wedge\Phi_1^k+\Phi_1^j\wedge\Phi_2^k+\Phi_2^j
  \wedge\Phi_1^k\bigg\}\,.
\end{array}\end{equation}
We expand
$\delta_{G\times
G}\{\Phi_2^j\wedge\Phi_2^k\}=\sum_i\{D{}_{1,ijk}\Phi_1^i+D{}_{2,ijk}\Phi_2^i\}$. Then,
taking into
account the normalizing factor of $2$ in Equation (\ref{eqn-1.a}) which dually  yields a
factor of $\frac12$ on the
inner product for $\Lambda^1(G\times G)$ and a factor of $\frac14$ on the inner product for
$\Lambda^2(G\times G)$, one has:
\begin{eqnarray*}
&&{\textstyle\frac12}D{}_{1,ijk}\text{vol}(G\times G)
=(\delta_{G\times G}\{\Phi_2^j\wedge\Phi_2^k\},\Phi_1^i)_{L^2(\Lambda^1(G\times G))}\\
&&\quad=(\Phi_2^j\wedge\Phi_2^k,d_{G\times G}\{\Phi_1^i\})_{L^2(\Lambda^2(G\times G))}=0,\\
&&{\textstyle\frac12}D{}_{2,ijk}\text{vol}(G\times G)
=(\delta_{G\times G}\{\Phi_2^j\wedge\Phi_2^k\},\Phi_2^i)_{L^2(\Lambda^1(G\times G))}\\
&&\quad=(\Phi_2^j\wedge\Phi_2^k,d_{G\times G}\{\Phi_2^i\})_{L^2(\Lambda^2(G\times G))}
={\textstyle\frac14}C_{ijk}\text{vol}(G\times G)\,.
\end{eqnarray*}
This shows that
\begin{equation}\label{eqn-4.g}
D{}_{1,ijk}=0\quad\text{and}\quad D{}_{2,ijk}=\textstyle\frac12C_{ijk}\,.
\end{equation}
Equations (\ref{eqn-4.d}), (\ref{eqn-4.e}),
(\ref{eqn-4.f}), and (\ref{eqn-4.g}) yield:
\begin{eqnarray*}
&&\Delta_{G\times G}(\Phi_2^i)=\delta_{G\times G}d_{G\times G}\{\Phi_2^i\}
={\textstyle\frac12}\sum_{l,j<k}C_{ljk}C_{ijk}\Phi_2^l=\textstyle\frac12\lambda_i\Phi_2^i\,.
\end{eqnarray*}
Similarly
\begin{eqnarray*}
&&\delta_{G\times G}\{\Phi_1^j\wedge\Phi_1^k\}=\delta_{G\times G}\{\Phi_2^j\wedge\Phi_1^k\}
=\delta_{G\times G}\{\Phi_1^j\wedge\Phi_2^k\}=\textstyle\frac12\sum_lC_{ljk}\Phi_1^l
\end{eqnarray*}
and thus $\Delta_{G\times G}^1\{\Phi_1^i\}=\frac32\lambda_i\Phi_1^i$.
\end{proof}

\section{Eigen forms whose pull-back has many non-zero Fourier coefficients}\label{sect-5}

Let $S^3$ be the unit sphere in the quaternions $\mathbb{H}=\mathbb{R}^4$; this is a
compact connected Lie group
and the standard round metric is the only bi-invariant metric on
$S^3$ modulo rescaling.  Fix
$$0\ne f\in E_{\lambda_0}(\Delta_{S^3}^0)$$ with
$\lambda_0\ne0$. Since the first cohomology group of $S^3$ is
trivial, there are no non-trivial harmonic $1$-forms on $S^3$.
Thus we may choose $$0\ne\phi\in\Lambda_L^1(S^3)\cap
E_{\lambda_1}(\Delta_{S^3}^1)$$ for some $\lambda_1>0$; we
refer to \cite{GLP96,LC05} for additional details concerning the
spectral geometry of $S^3$; $S^3$ could be replaced by any
non-Abelian compact connected Lie group in this construction.

We first prove Theorem \ref{thm-1.4} in the special case that
$p=1$.  Suppose that $\mu_0=2k$. Choose real numbers
$0<t_1<...<t_k<1$. Choose
$s_1>...>s_k>1$ so $$s_\alpha \lambda_0+t_\alpha
\lambda_1=\lambda_0+\lambda_1\quad\text{for}\quad1\le\alpha\le k\,.$$
Let
$G_\alpha $ be $S^3$ with the rescaled metric $ds_{G_\alpha
}^2:=t_\alpha ^{-1}ds_{S^3}^2$ and let $\phi^\alpha
=\phi\in\Lambda_L^1(G_\alpha )$. Let $\bar G_\alpha $ be $S^3$
with the rescaled metric $ds_{\bar G_\alpha }^2:=s_\alpha
^{-1}ds_{S^3}^2$ and let $f_\alpha =f\in C^\infty(\bar G_\alpha
)$. After taking into account the effect of the rescaling, we have
$$f_\alpha \in E_{s_\alpha \lambda_0}(\Delta_{\bar G_\alpha }^0),\quad
  df_\alpha \in E_{s_\alpha \lambda_0}(\Delta_{\bar G_\alpha
}^1),\quad\text{and}\quad\phi^\alpha \in E_{t_\alpha \lambda_1}(\Delta_{G_\alpha }^1)\,.$$
Let $G=G_1\times...\times G_k\times \bar G_1\times...\times \bar G_k$. Decompose
$m^*(\phi^\alpha )=\Phi^\alpha _1+\Phi^\alpha _2$.
Let $\psi:=\sum_\alpha f_\alpha\phi^\alpha$. As the structures decouple, one has:
$$\Delta_G^1\{\psi\}=\sum_\alpha (s_\alpha \lambda_0+t_\alpha\lambda_1)f_\alpha
\phi^\alpha=(\lambda_0+\lambda_1)\psi \,.$$ We can apply Theorem \ref{thm-1.3} to see
\begin{eqnarray*}
&&\Delta_{G\times G}^1m^*\psi=\sum_\alpha \{(s_\alpha \lambda_0+\textstyle\frac32t_\alpha
\lambda_1)m^*f_\alpha
\cdot\Phi_1^\alpha
  +(s_\alpha \lambda_0+\textstyle\frac12t_\alpha \lambda_1)m^*f_\alpha \cdot\Phi_2^\alpha
\}\\
&&\qquad=\sum_\alpha \{(\lambda_0+\lambda_1+
\textstyle\frac12t_\alpha \lambda_1)m^*f_\alpha \cdot\Phi_1^\alpha
  +(\lambda_0+\lambda_1-\textstyle\frac12t_\alpha \lambda_1)m^*f_\alpha
\cdot\Phi_2^\alpha \}\,.
\end{eqnarray*}
The computations performed above then yield $\psi\in E_{\lambda_0+\lambda_1}(\Delta_G^1)$.
Furthermore:
\begin{eqnarray*}
&&m^*(f_\alpha )\Phi_1^\alpha \in E_{\lambda_0+\lambda_1+\frac12t_\alpha
\lambda_1}(\Delta^1_{{G\times G}}),\\
&&m^*(f_\alpha )\Phi_2^\alpha \in
E_{\lambda_0+\lambda_1-\frac12t_\alpha
\lambda_1}(\Delta^1_{{G\times G}})\,.
\end{eqnarray*}
Since $0<t_1<...<t_k$, $m^*\psi$ has a Fourier decomposition which
involves $2k=\mu_0$ distinct eigenvalues. This establishes Theorem
\ref{thm-1.4} if $p=1$ and if $\mu_0$ is even.

If $\mu_0=2k+1$ is odd, we choose $s_0$ so
$s_0\lambda_0=\lambda_0+\lambda_1$. Then $f_0\in
E_{\lambda_1+\lambda_2}(\Delta_{\bar G_0}^0)$. We apply the
construction described above to $G=G_1\times...\times G_k\times
\bar G_0\times...\times \bar G_k$ and to
$\psi=df_0+f_1\phi^1+...+f_k\phi^k$; the latter factors are not
present if $\mu_0=1$. Since $m^*df_0\in
E_{\lambda_0+\lambda_1}(\Delta_{G\times G}^1)$, there are $2k+1$
distinct eigenvalues which are involved in the Fourier
decomposition of $\psi$. This completes the proof of Theorem
\ref{thm-1.4} if $p=1$. We take the product of $G$ with circles
$S^1$ and replace $\phi$ by $\phi\wedge d\theta_1\wedge...\wedge
d\theta_p$, where $\theta_\beta$ is the usual periodic parameter
on $S^1$, to complete the proof if $p\ge1$.\hfill\qedbox

\section{Harmonic forms}\label{sect-1.6}

Before beginning the proof of Theorem \ref{thm-1.5}, we must establish some technical results. Let
$\gg_L$ be the Lie algebra of left invariant vector fields on $G$. The following results are well
known; we sketch the proofs briefly:

\begin{lemma}\label{lem-6.2}
Let $ds^2_G$ be a bi-invariant metric on a compact connected Lie
group $G$.
\begin{enumerate}
\item If $\theta\in E_0(\Delta_G^n)$, then $\theta$ is bi-invariant.
\item If $\eta$ is a bi-invariant vector field, then $\nabla\eta=0$.
\item Let $\theta\in\Lambda_L^1(G)$. If $d\theta=0$, then $\nabla\theta=0$.
\item If $\Theta\in\Lambda^n(E_0(\Delta_G^1))$, then $\nabla\Theta=0$ and $\Theta\in E_0(\Delta_G^n)$.
\end{enumerate}
\end{lemma}

\begin{proof} The Hodge-DeRham theorem provides a natural identification of
$E_0(\Delta_G^n)$ with the cohomology group $H^n(G;\mathbb{C})$. In particular, this identification
is compatible with the action of $L_g^*$ and $R_g^*$. Since $G$ is connected, $L_g^*$ and
$R_g^*$ act trivially on $H^n(G;\mathbb{C})$ and hence on $E_0(\Delta_G^n)$. Assertion (1) follows.

To prove Assertion (2), we use well known facts concerning bi-invariant metrics on Lie groups; see, for example,
\cite{L97}.  Let
$\exp(t\xi)$ be the integral curve through the identity for $\xi\in\gg_L(G)$.  Let $\eta$ be bi-invariant. Assertion
(2) follows as:
$$\nabla_\xi\eta=\textstyle\frac12[\xi,\eta]
   =\textstyle\frac12\partial_t\left\{(L_{\exp(t\xi)})_*(R_{\exp(-t\xi)})_*\eta\right\}|_{t=0}
=\partial_t\eta|_{t=0}=0\,.
$$

Let $\theta\in\Lambda_L^1(G)$ with $d\theta=0$. Since $\delta\theta$ is left
invariant, $\delta\theta=c$ is constant. Since $\Delta_G^0c=0$, $\delta\theta=0$. Thus
$\theta$ is harmonic and hence bi-invariant. We use the metric to raise and lower indices and
identify the tangent and cotangent spaces. Let $\eta$ be the corresponding dual bi-invariant
vector field. By Assertion (2), $\eta$ is parallel. Thus, dually, $\theta$ is parallel. This
proves Assertion (3).

Let $\Theta\in\Lambda^n(E_0(\Delta_G^1))$. Then there are constants $a_I$ and harmonic $1$-forms $\theta^i_L$ so
$$\Theta=\sum_{|I|=n}a_I\theta^{i_1}_L\wedge...\wedge\theta^{i_n}_L\,.$$
By assertion (3), $\nabla\theta^i_L=0$. Consequently $\nabla\Theta=0$. On the other hand, one has
$$d+\delta=\sum_i\{\operatorname{ext}(e^i)-\operatorname{int}(e^i)\}\nabla_{e_i}$$
where $\{e_i\}$ and $\{e^i\}$ are dual orthonormal frames for $TG$ and $T^*G$ and where
$\operatorname{ext}(\cdot)$ and
$\operatorname{int}(\cdot)$ denote exterior and interior multiplication.  Thus parallel forms are necessarily
harmonic. Assertion (4) follows.
\end{proof}

We distinguish the two factors in the product to decompose
$$\Lambda^n(G\times G)=\oplus_{p+q=n}\Lambda^p(G_1)\otimes\Lambda^q(G_2)\,.$$
We let $\pi_{p,q}$
denote orthogonal projection on the various components. The K\"unneth formula shows
$$H^n(G\times G;\mathbb{C})=\oplus_{p+q=n}H^p(G_1;\mathbb{C})\otimes H^q(G_2;\mathbb{C})$$
and, as we have taken a product metric on $G\times G$, we have a corresponding decomposition in the
geometric context:
\begin{eqnarray*}
&&C^\infty(\Lambda^n(G\times G))=\oplus_{p+q=n}C^\infty\{\Lambda^p(G_1)\otimes\Lambda^q(G_2)\},\\
&&\Delta_{G\times G}^n=\oplus_{p+q=n}\{\Delta_{G_1}^p\otimes\id+\id\otimes\Delta_{G_2}^q\},\\
&&E_0(\Delta_{G\times G}^n)=\oplus_{p+q=n}\{E_0(\Delta_{G_1}^p)\otimes E_0(\Delta_{G_2}^q)\}\,.
\end{eqnarray*}

Assertion (1) of Theorem \ref{thm-1.5} follows from Lemma \ref{lem-6.2}. To
prove Assertion (2) of Theorem \ref{thm-1.5}, suppose
$$\phi=\sum_{|I|=n}a_I\theta^{i_1}\wedge...\wedge\theta^{i_n}\in\Lambda^n(E_0(\Delta_G^1))\text{ for }
\theta^j\in E_0(\Delta^1_G)\,.$$
As $\theta^j$ is bi-invariant,
$\theta^j\oplus\theta^j\in\tilde\Lambda^1(G\times G)$. Since
$m^*\theta(1,1)=\theta^j(1,1)\oplus\theta^j(1,1)$, $m^*\theta^j=\theta^j\oplus\theta^j$.
As $d\theta^j=0$,
$dm^*\theta^j=m^*d\theta^j=0$. Thus $m^*\theta^j\in E_0(\Delta_{G\times G}^1)$ so
$m^*\phi\in\Lambda^n(E_0(\Delta_{G\times G}^1))$ is harmonic.

Conversely suppose that $m^*\phi\in E_0(\Delta_{G\times G}^n)$. We then have $\pi_{0,n}m^*\phi$ is harmonic.
Since $\pi_{0,n}m^*\phi=\sigma_2^*\phi=0\oplus\phi$, $\phi$ is harmonic and hence bi-invariant. Decompose
$\phi=\sum_{|I|=n}a_I\phi_L^I$ as a sum of left invariant $n$-forms where the coefficients $a_I$ are constant.
As $m^*\phi$ is harmonic, $m^*\phi$ is left invariant and decomposes in the form:
$$m^*\phi=\sum_{0<i_1<...<i_n<\dim(G)}a_{i_1...i_n}(\phi_L^{i_1}\oplus0+0\oplus\phi_L^{i_1})
\wedge...\wedge(\phi_L^{i_n}\oplus0+0\oplus\phi_L^{i_n})\,.$$
Choose the indexing convention so $\{\phi^1,...,\phi^k\}$ is an orthonormal basis for
$E_0(\Delta_G^1)$ and so $\{\phi^{k+1},...,\phi^{\dim(G)}\}$ completes the set to an orthonormal basis for
$\Lambda_L^1(G)$. We suppose that
$\phi\notin\Lambda^n(\phi^1,...,\phi^k)$ and argue for a contradiction. Choose $a$ minimal so
$a_{i_1,...,i_a,j_1,...,j_b}\ne0$ where $i_a\le k$ and
$k<j_1<j_2<...<j_b$. By hypothesis $a<n$ so $b\ge1$. Let
\begin{eqnarray*}
&&\phi_0:=\phi_L^{i_1}\wedge...\wedge\phi_L^{i_a}\wedge\phi_L^{j_2}\wedge...\wedge\phi_L^{j_b},\\
&&\phi:=\tilde\phi\wedge\phi_0+\text{other terms}
\end{eqnarray*}
where the other terms do not involve the monomial $\phi_0$ and where $0\ne\tilde\phi\notin
E_0(\Delta_G^1)$; $\tilde\phi=\operatorname{int}(\phi_0)\phi\in\Lambda_L^1(G)$. We may then expand
\begin{eqnarray*}
&&\pi_{1,n-1}m^*\phi=(\tilde\phi\oplus0)\wedge(0\oplus\phi_0)+\text{other terms},\\
&&d\pi_{1,n-1}m^*\phi=(d\tilde\phi\oplus0)\wedge(0\oplus\phi_0)+\text{other terms}\,.
\end{eqnarray*}
Consequently $d\tilde\phi=0$ since this is the only term of bi-degree $(2,n-1)$
multiplied by $0\oplus\phi_0$; one then has $d\tilde\phi\oplus0=\operatorname{int}(0\oplus\phi_0)dm^*\phi$. By Lemma
\ref{lem-6.2},
$\tilde\phi\in E_0(\Delta_G^1)$. The contradiction completes the proof of Assertion (2).
Assertion (3) is an immediate consequence of Assertion (2) since $E_0(\Delta_G^1)=\{0\}$ if $G$ is
simply connected.\hfill\qedbox

\section{Finite Fourier series for general left invariant metrics}\label{sect-1.7}

Let $ds^2_G$ be a left invariant metric on a compact Lie group $G$ and let $ds^2_{G\times G}$ be a
left invariant metric on $G\times G$. We impose no relationship between the two metrics; in particular, we do
not assume that the multiplication map $m$ is a Riemannian submersion any more; it is an interesting question in
its own right when this is possible and we shall investigate this question in more detail in a subsequent paper.

We begin the proof of Theorem \ref{thm-1.6} with a technical result:
\begin{lemma}\label{lem-7.1}
Let $G$ be a compact Lie group. Let $H$ be left invariant subspace of
$L^2(\Lambda^p(G))$. Then there is a bi-invariant subspace
$H_1\subset L^2(\Lambda^p(G))$ which contains $H$ so that
$\dim(H_1)\le{\textstyle\binom{\dim\{G\}}{p}}\dim(H)^2$.
\end{lemma}

\begin{proof} The Lemma is immediate if $\dim(H)=\infty$ so we may suppose $H$ is finite
dimensional. By decomposing
$H=\oplus_in_iV_{\rho_i}$ into irreducible representations,  we may assume without
loss of generality that $H=V_\rho$ where $V_\rho$ is an irreducible left representation space
for $G$ in the
proof of Lemma \ref{lem-7.1}.  We apply the Peter-Weyl theorem and use Equation (\ref{eqn-3.a}).
Let
$$H_\rho^p=\oplus_{|I|=p}H_\rho\cdot\Phi_L^I={\textstyle\binom{\dim\{G\}}{p}}\dim(V_\rho)V_\rho\,.$$
Then $H\subset H_\rho^p$. Since
left and right multiplication commute, $H_\rho^p\cdot R_g$ is isomorphic to $H_\rho^p$ for any
$g\in G$. Since $H_\rho^p$ contains all representations isomorphic to $V_\rho$, $H_\rho^p\cdot R_g=H_\rho^p$ is
right invariant as well.
\end{proof}

Let
$\phi\in E_\lambda(\Delta_G^p)$ be an eigen
$p$-form. Apply Lemma \ref{lem-7.1} to choose a subspace $H_1\subset L^2(\Lambda^p(G))$ which is
left and right invariant under the action of
$G$, which contains
$E_\lambda(\Delta_G^p)$, and which satisfies:
$$\dim(H_1)\le{\textstyle\binom{\dim\{G\}}{p}}\dim\{E_\lambda(\Delta_G^p)\}^2\,.$$
Let $\tilde m(a,b)=ab^{-1}$. Then $\tilde m L_{g_1,g_2}=L_{g_1}R_{g_2^{-1}}\tilde m$. Consequently
$\tilde m^*H_1$ is a finite dimensional subspace of $L^2(\Lambda^p(G\times G))$ which invariant
under left multiplication in the group. Apply Lemma \ref{lem-7.1} to choose a subspace $H_2\subset
L^2(G\times G)$ which is left and right under the action of $G\times G$, which contains
$\tilde m^*H_1$, and which satisfies:
$$\dim(H_2)\le{\textstyle\binom{2\dim\{G\}}{p}}{\textstyle\binom{\dim\{G\}}{p}}^2\dim\{E_\lambda(\Delta_G^p)\}^4
\,.$$
Set $\psi(x,y)=(x,y^{-1})$. Then
$\psi^*H_2$ is still bi-invariant and in particular is left invariant. Since
$m^*\phi\in\psi^*H_2$, Lemma
\ref{lem-3.1} can now be applied to show that
$$\mu(m^*\phi)\le
{\textstyle\binom{2\dim\{G\}}{p}}^2{\textstyle\binom{\dim\{G\}}{p}}^2\dim\{E_\lambda(\Delta_G^p)\}^4\,.
$$ Theorem
\ref{thm-1.6} now follows.\hfill\qedbox

\section*{Acknowledgments}
Research of C. Dunn partially supported by a CSUSB faculty research
grant. Research of P. Gilkey partially supported by the Max Planck
Institute in the Mathematical Sciences (Leipzig, Germany) and the
Program on Spectral Theory and Partial Differential Equations of the
Newton Institute (Cambridge UK). Research of J.H. Park partially
supported by the Korea Research Foundation Grant funded by the
Korean Government (MOEHRD) KRF-2005-204-C00007.

\end{document}